\documentclass{amsart}

\usepackage{bm,url,amsmath,amssymb,amsthm}

\newcommand{\Q}{\mathbb{Q}}
\newcommand{\C}{\mathbb{C}}
\newcommand{\R}{\mathbb{R}}
\newcommand{\Z}{\mathbb{Z}}
\newcommand{\h}{\mathbb{H}}
\newcommand{\w}{\mathcal{H}}
\newcommand{\s}{\Lambda}
\newcommand{\alg}{\mathcal{A}}
\newcommand{\ord}{\mathcal{O}}

\newtheorem{theorem}{Theorem}[section]
\newtheorem{corollary}[theorem]{Corollary}
\newtheorem{proposition}[theorem]{Proposition}
\newtheorem{lemma}[theorem]{Lemma}

\theoremstyle{definition}

\theoremstyle{remark}

\numberwithin{equation}{section}

\begin{document}

\title{A Mordell Inequality for Lattices \linebreak over Maximal Orders}


\author{Stephanie Vance}
\address{\textit{Current Address:} School of Sciences\\
Adams State College\\
208 Edgemont Blvd.\\
Alamosa, CO 81102}
\address{\textit{Former Address:} Department of Mathematics\\ University of Washington\\ Box 354350\\ Seattle, WA
98195}
\email{slvance@adams.edu}
\thanks{The author was supported by an ARCS Foundation fellowship and a research assistantship funded by Microsoft Research}


\date{}

\dedicatory{}

\begin{abstract}
In this paper we prove an analogue of Mordell's inequality for lattices in finite-dimensional complex or quaternionic Hermitian space that are modules over a maximal order in an imaginary quadratic number field or a totally definite rational quaternion algebra.  This inequality implies that the $16$-dimensional Barnes-Wall lattice has optimal density among all $16$-dimensional lattices with Hurwitz structures.  
\end{abstract}

\maketitle

\section{Introduction}

Lattice sphere packings in $n$-dimensional Euclidean space are configurations of congruent spheres with disjoint interiors in which the sphere centers form a lattice.  These geometric objects and more specifically their densities, i.e., the proportions of space covered by the spheres, have been an active area of research over the past several centuries.  One particular problem in this research area that has received much attention is the \it $n$-dimensional lattice sphere packing problem\normalfont.  For this problem one must determine the optimal density of a lattice in $n$-dimensional Euclidean space, with the density of a lattice defined to be the density of the lattice sphere packing obtained by centering at each lattice point a sphere with radius equal to half the length of the shortest non-zero lattice vectors.  Remarkably this problem remains open for all dimensions $n>8$ with the exception of dimension $24$; see \cite{CS1} and \cite{CK} for more details regarding the known optimal lattices in dimensions $n\leq 8$ and dimension $24$. 

In most low dimensions divisible by $2$ or $4$, at least one of the densest known lattices has an Eisenstein or a Hurwitz structure respectively; i.e., at least one of these lattices is isometric to a lattice in a complex or quaternionic Hermitian space that is a module over the Eisenstein integers $\mathcal{E} = \Z[\frac{1+\sqrt{-3}}{2}]$ in the complex field $\C$ or the Hurwitz integers $\w = \Z[i,j,\frac{1+i+j+k}{2}]$ in the quaternion skew-field
$$\h = \lbrace a+bi+cj+dk : a,b,c,d\in\R,\ i^2=j^2=-1 \text{ and } ij = -ji = k\rbrace,$$ respectively.  
For example, the densest known lattices in every even dimension up to $24$ all have an Eisenstein structure, and the densest known lattices in dimensions $4$, $8$, $16$, and $24$ all have a Hurwitz structure (note that the densest known $12$-dimensional lattice does not have a Hurwitz structure).
The existence of this extra algebraic structure among the densest known lattices makes it natural to consider the $2m$-dimensional Eisenstein and $4m$-dimensional Hurwitz lattice sphere packing problems, i.e., determine the optimal density of $2m$ and $4m$-dimensional Eisenstein and Hurwitz lattices. 
Even if the densest lattices in certain dimensions do not have an Eisenstein or a Hurwitz structure, it would still be interesting to determine the optimal density of the lattices that do satisfy these additional algebraic constraints.

In this paper we address the $2m$ and $4m$-dimensional Eisenstein and Hurwitz lattice sphere packing problems simultaneously by considering Eisenstein and Hurwitz lattices in the more general context of $\ord$-lattices, i.e., lattices in finite-dimensional complex or quaternionic Hermitian space that are modules over a maximal $\Z$-order $\ord$ in an imaginary quadratic number field or a totally definite quaternion $\Q$-algebra, respectively.    
In the next section we prove several propositions for $\ord$-lattices and their determinants, i.e., the squared volumes of their fundamental regions.   
These propositions are used in Section 3 to prove a Mordell inequality for $\ord$-lattices which relates the optimal value of the Hermite invariant of $\ord$-lattices in two consecutive complex or quaternionic dimensions. Note that the Hermite invariant of an $n$-dimensional lattice $\s$ is directly proportional to the $(n/2)^\mathrm{th}$ power of its density and is given by the formula 
$$\gamma(\s) = \frac{\mathrm{N}(\s)}{{\det(\s)}^{1/n}},$$
where $\mathrm{N}(\s)$ denotes the norm of $\s$ (i.e., the norm of the shortest non-zero vectors in $\s$) and $\det(\s)$ denotes the determinant of $\s$.   
Then in Section 4 we use the Eisenstein and Hurwitz versions of this Mordell inequality to obtain upper bounds for the Hermite invariants of Eisenstein and Hurwitz lattices in low dimensions, and we show that the $16$-dimensional Barnes-Wall lattice has optimal density as a $16$-dimensional lattice with a Hurwitz structure.  

 While the focus of this paper is on Eisenstein and Hurwitz lattices, we also consider in Section 5 Gaussian lattices in $\C^m$ (i.e., $\ord = \Z[i]$) and the lattices in $\h^m$ that are modules over the ring $\ord = \Z[1,i,\frac{1+j}{2},\frac{i+k}{2}]$.  As with the Eisenstein and Hurwitz cases discussed above, some of the densest known  lattices in dimensions $2m$ and $4m$ have one of these $\ord$-lattice structures.   Note that all four types of $\ord$-lattices discussed in this paper can be regarded as $G$-lattices in the sense of \cite[ch.~13]{Ma}.  


\section{$\ord$-Lattices in Complex and Quaternionic Hermitian space}

The following notation is used throughout this paper.

\begin{itemize} 
  \item Let $\mathbb{K}$ denote either the complex field $\C$ or the quaternionic skew-field $\h$, and let $r$ denote the rank of $\mathbb{K}$ as an $\R$-algebra, i.e., $r = 2$ or $r = 4$.  Then let $x\mapsto\overline{x}$ denote complex conjugation when $\mathbb{K} = \C$ and quaternionic conjugation, i.e., $a+bi+cj+dk\mapsto a-bi-cj-dk$ for $a,b,c,d\in\R$, when $\mathbb{K} = \h$.
  \item When $\mathbb{K} = \C$ let $\alg$ denote an imaginary quadratic number field and when $\mathbb{K} = \h$ let $\alg$ denote a totally definite quaternion $\Q$-algebra. In both cases we shall identify $\mathbb{K}$ with $\R\otimes_\Q\alg$.
  \item Let $\ord$ denote the image in $\mathbb{K}$ of a maximal $\Z$-order in the $\Q$-algebra $\alg$; i.e., $\ord$ is a free $\Z$-submodule and subring of $\alg$ (where $\alg$ is identified with its image in $\mathbb{K}$) satisfying $\mathrm{rank}_\Z\ord = \mathrm{rank}_\R K$ and $\Q\ord = \alg$. Then let $D_{\ord} =  \det\left(\left(\frac{1}{2}\left(\alpha_i\overline{\alpha_j} + \alpha_j\overline{\alpha_i}\right)\right)_{1\leq i,j\leq r}\right)$, where $\alpha_1,\dots,\alpha_r$ denotes a $\Z$-basis for $\ord$. Note that $D_\ord$ is independent of the $\Z$-basis for $\ord$ used to compute it, and hence is an invariant of $\ord$. 
  \item Finally let $E_{rm}$ denote a left $m$-dimensional $\mathbb{K}$-vector space with a non-degenerate Hermitian product $h:\mathbb{K}\times\mathbb{K}\to\mathbb{K}$ (i.e., $h$ is left linear in the first variable and right conjugate linear in the second variable) and define the inner product of each pair of vectors $x,y\in E_{rm}$ to be
$$\langle x,y\rangle = \frac{1}{r}\mathrm{Tr}_{\mathbb{K}/\R}(h(x,y)) = \frac{1}{2}\left(h(x,y) + \overline{h(x,y)}\right)$$
so that $\mathrm{N}(x) = h(x,x) = \langle x,x\rangle$; i.e., the distances defined on $E_{rm}$ by $h(\cdot,\cdot)$ and  $\langle\cdot,\cdot\rangle$ are the same\footnote{The factor $\frac{1}{2}$ in the definition of the inner product $\langle\cdot,\cdot\rangle$ is  inserted to simplify later calculations and can be omitted without affecting the density of a lattice in $E_{rm}$.  One of the reasons for omitting this factor would be to ensure that every lattice $\s$ having the property that $h(x,y)\in\ord$ for every $x,y\in\s$ is necessarily integral, i.e., $\langle x,y\rangle\in\Z$ for every $x,y\in\s$.}.
\end{itemize}  
 
Using this notation we define an $rm$-dimensional $\ord$-lattice to be a left $\ord$-invariant lattice in $E_{rm}$, i.e., a lattice in $E_{rm}$ that is a left $\ord$-module.  In this section we prove several propositions concerning an $rm$-dimensional $\ord$-lattice $\s$ and its $\ord$-dual, which is defined to be the $\ord$-module 
$$\s^\# = \lbrace x\in E_{rm}: h(x,\s)\subseteq\ord\rbrace.$$

\begin{proposition}\label{prop: new} If $\s$ is an $\ord$-lattice in $E_{rm}$, then $\s^\#$ is also an $\ord$-lattice in $E_{rm}$, and the determinants of these two lattices satisfy the product formula
$$\det(\s)\det(\s^\#) = D_{\ord}^{2m}.$$
\end{proposition}

To prove Proposition \ref{prop: new} we use the following lemma to deal with the case when $\s$ is a non-free $\ord$-lattice (i.e., when $\s$ is not a free $\ord$-module).  


\begin{lemma}\label{lemma: basic}
For every $\ord$-lattice $\s$ in $E_{rm}$ there exists two free $\ord$-lattices $\s_1$ and $\s_2$ in $E_{rm}$ such that $\s_1\subseteq\s\subseteq\s_2$.   
\end{lemma}

\begin{proof}An $\ord$-lattice $\s$ generates $E_{rm}$ as a real vector space and so must contain a $\mathbb{K}$-vector space basis for $E_{rm}$.  Hence we can choose $\s_1$ to be the free $\ord$-lattice generated over $\ord$ by such a $\mathbb{K}$-vector space basis contained in $\s$.  Now consider the two $\Q$-vector subspaces $\Q\s_1$ and $\Q\s$ of $E_{rm}$.  Both subspaces have dimension $rm$, and since $\Q\s_1\subseteq\Q\s$ they must be equal.  Thus $\s$ is a finitely generated submodule of $\Q\s_1$, and so we can choose an $a\in\Z$ such that $\s\subseteq a^{-1}\s_1$, i.e., $\s$ is contained in the free $\ord$-lattice $\s_2 = a^{-1}\s_1$.

\end{proof}

\begin{proof}[Proof of Proposition \ref{prop: new}]
We first verify that the $\ord$-dual of every $\ord$-lattice in $E_{rm}$ is also an $\ord$-lattice in $E_{rm}$.  
Every $\ord$-basis for a free $\ord$-lattice in $E_{rm}$ is necessarily a $\mathbb{K}$-vector space basis for $E_{rm}$, say $v_1,\dots,v_m$, 
 and the dual $\mathbb{K}$-vector space basis $v_1^\#,\dots,v_m^\#$ satisfying $h(v_i^\#,v_j) = \delta_{i,j}$ is an $\ord$-basis for the $\ord$-dual.  Hence the $\ord$-dual of a free $\ord$-lattice in $E_{rm}$ is also a free $\ord$-lattice in $E_{rm}$.  In particular, if $\s_1$ and $\s_2$ are the two free $\ord$-lattices in $E_{rm}$ given by Lemma \ref{lemma: basic}, then both $\s_1^\#$ and $\s_2^\#$ are free $\ord$-lattices in $E_{rm}$.  The inclusion $\s_1\subseteq\s\subseteq\s_2$ then implies that $\s_2^\#\subseteq\s^\#\subseteq\s_1^\#$, and so $\s^\#$ must also be an $\ord$-lattice in $E_{rm}$.

To prove the determinant identity we first consider the case when $\s$ is a free $\ord$-lattice with $\ord$-basis $v_1,\dots,v_m$.  Letting $\alpha_1,\dots,\alpha_r$ denote a $\Z$-basis for $\ord$, the vectors $\lbrace \alpha_i v_j\rbrace_{1\leq i\leq r,1\leq j\leq m}$ form a $\Z$-basis for $\s$ and the vectors $\lbrace \alpha_i v_j^\#\rbrace_{1\leq i\leq r,1\leq j\leq m}$ form a $\Z$-basis for $\s^{\#}$.
Write the vectors in these two $\Z$-bases using coordinates with respect to a fixed orthonormal basis for $E_{rm}$ (the basis being orthonormal with respect to the inner product $\langle\cdot,\cdot\rangle$), and let $M$ be an $rm\times rm$ matrix such that for $1\leq i\leq r$ and $1\leq j\leq m$ the $(r(j-1)+i)^{\mathrm{th}}$ row of $M$ is the vector $\alpha_i v_j$.  Similarly let $N$ be an $rm\times rm$ matrix such that for $1\leq i\leq r$ and $1\leq j\leq m$ the $(r(j-1)+i)^{\mathrm{th}}$ row of $N$ is the vector $\alpha_i v_j^\#$.  Observe that the two matrices $M$ and $N$ satisfy $\det(\s) = \det(MM^\mathrm{T})$ and $\det(\s^\#) = \det(NN^\mathrm{T})$.  Hence the discriminant ideals\footnote{The discriminant ideal of an $n$-dimensional lattice $\s$ is equal to the ideal in $\Z$ consisting of the determinants of all $n\times n$ Gram matrices $\left(\langle w_i,w_j\rangle\right)$ where $w_1,\dots,w_n\in\s$.  Note that since $\s$ is a free $\Z$-module of rank $n$, it is not hard to show that $d(\s) = \Z\det(\s)$.}    
$d(\s)$ and $d(\s^\#)$ satisfy 
\begin{eqnarray*}
d\left(\s\right)d\left(\s^{\#}\right)& = &\Z\det(MM^\mathrm{T})\det(NN^\mathrm{T})\\
&= & \Z\det(MN^{\mathrm{T}})^2.
\end{eqnarray*}
(Note that we are computing the product of discriminant ideals here rather than the lattice determinants, so that our calculations can be generalized for the non-free case.)

We can compute the determinant of the matrix $MN^{\mathrm{T}}$ using the fact that for $1\leq i_1,i_2\leq r$ and $1\leq j_1,j_2\leq m$, the $(r(j_1-1)+i_1,r(j_2-1)+ i_2)^{\mathrm{th}}$  entry in $MN^\mathrm{T}$ is equal to  
$$\langle\alpha_{i_1}v_{j_1},\alpha_{i_2}v_{j_2}^\#\rangle = \frac{1}{2}\left(h(\alpha_{i_1}v_{j_1},\alpha_{i_2}v_{j_2}^\#) +  h(\alpha_{i_2}v_{j_2}^\#,\alpha_{i_1}v_{j_1})\right)= \frac{1}{2}\left(\alpha_{i_1}\overline{\alpha_{i_2}}+\alpha_{i_2}\overline{\alpha_{i_1}}\right)\delta_{j_1,j_2}.$$  In particular, $MN^{\mathrm{T}}$ is an $rm\times rm$ block-diagonal matrix with $m$ blocks all equal to the $r\times r$ matrix 
$\left(\frac{1}{2}\left(\alpha_i\overline{\alpha_{j}}+\alpha_{j}\overline{\alpha_{i}}\right)\right)_{1\leq i,j\leq r}$ whose determinant is $D_\ord$.
Therefore, 
$$d(\s)d\left(\s^{\#}\right) = \Z\left(\prod_{s = 1}^mD_\ord\right)^2 =  \Z{D_\ord}^{2m},$$ and from these equalities the determinant identity follows.

Now suppose that $\s$ is a non-free $\ord$-lattice.  Notice that in the preceding definitions and calculations we may replace $\Z$  by any localization, in particular by $\Z_P$, where $P$ is a prime ideal in $\Z$.  Letting $S = \Z\backslash P$ we have  
$$(\s_P)^\# = \lbrace x\in K: \forall s\in S, h(x,s^{-1}\s)\subseteq S^{-1}\ord\rbrace,$$
and since $h(x,s^{-1}\s) = s^{-1}h(x,\s)$, we may replace $h(x,s^{-1}\s)$ by $h(x,\s)$ in the formula above, showing  that $(\s_P)^\# = (\s^\#)_P$.   
Moreover, since $\s_P$ is a free $\ord_P$-module (\cite{AG}, Proposition 3.7) generated by a $\mathbb{K}$-vector space basis of $E_{rm}$, the product of the discriminant ideals $d(\s_P)$ and $d((\s_P)^\#)$ is equal to 
$\Z_PD_{\ord_{P}}^{2m} = \Z_PD_\ord^{2m}$ (here $D_{\ord_P}$ is equal to the determinant of the Gram matrix corresponding to a $\Z_P$-basis for $\ord_P$).   Therefore for every prime ideal $P$ in $\Z$  we have $$\Z_Pd(\s)d(\s^\#) = d(\s_P)d((\s^\#)_P) = d(\s_P)d((\s_P)^\#) = \Z_PD_\ord^{2m},$$
implying that $d(\s)d(\s^\#) = \Z D_\ord^{2m}$.  The determinant identity $\det(\s)\det(\s^\#) = D_\ord^{2m}$ follows from this last equality.


\end{proof}

Observe from the definition of the $\ord$-dual of an $\ord$-lattice $\s$ that we always have the containment $\s\subseteq\s^{\#\#}$.  Now since the determinant identity in Proposition \ref{prop: new} implies that $\det(\s) = \det(\s^{\#\#})$, it then follows that $\s = \s^{\#\#}$.

We conclude this section with another determinant identity for $\ord$-lattices that involves the intersection of an $rm$-dimensional $\ord$-lattice $\s$ with a $\mathbb{K}$-vector subspace $F$ of $E_{rm}$ and the projection of $\s$ onto $F^{\perp}$, i.e., the $\mathbb{K}$-vector subspace of $E_{rm}$ perpendicular to $F$ with respect to $\langle\cdot,\cdot\rangle$ (or equivalently $h(\cdot,\cdot)$ due to the relationship between $\langle\cdot,\cdot\rangle$ and $h(\cdot,\cdot)$).  Note that below we refer to $\s\cap F$ as a \it relative lattice \normalfont in $F$ because it is a lattice in a $\mathbb{K}$-vector subspace of $F$ that may not be full-dimensional.  
Also, we regard both $\mathbb{K}$-vector subspaces $F$ and $F^\perp$ of $E_{rm}$ as $\mathbb{K}$-Hermitian spaces (and hence call them $\mathbb{K}$-Hermitian subspaces) because we can restrict the Hermitian product $h(\cdot,\cdot)$ to $F\times F$ and $F^\perp\times F^\perp$, respectively.

\begin{lemma} \label{prop: reldet}
Let $\s$ be an $rm$-dimensional $\ord$-lattice and let $F$ be a $\mathbb{K}$-Hermitian subspace of $E_{rm}$.  
\begin{enumerate}
  \item The relative $\ord$-lattice $\s\cap F$ is an $\ord$-lattice in $F$ if and only if $\pi_{F^\perp}(\s)$ is an $\ord$-lattice in $F^\perp$.  
  \item If $\s \cap F$ is an $\ord$-lattice in $F$, then $\det(\s) = \det(\s\cap{F})\det(\pi_{F^{\perp}}(\s)).$
\end{enumerate}  
\end{lemma}
\begin{proof}
This lemma is proved as Proposition 1.2.9 in \cite{Ma} for lattices in finite-dimensional Euclidean space.  Because $\s\cap F$ and $\pi_{F^\perp}(\s)$ are both $\ord$-modules, the $\ord$-lattice version readily follows.
\end{proof}

\begin{proposition}\label{prop: intersection}
If $\s$ is an $rm$-dimensional $\ord$-lattice and $F$ is a $\mathbb{K}$-Hermitian subspace in $E_{rm}$, then $\s\cap F$ is an $\ord$-lattice in $F$ if and only if $\s^\#\cap F^\perp$ is an $\ord$-lattice in $F^\perp$.  Moreover if these conditions hold and $s = \mathrm{dim}_K F$, then
$$\det(\s) = \det(\s\cap{F})\det(\s^\#\cap F^\perp)^{-1} {D_\ord}^{2(m-s)}.$$
\end{proposition}

\begin{proof}  
For every $x\in F^\perp$ and $y\in \s$ we can write $h(x,y) = h(x,\pi_{F}(y) + \pi_{F^\perp}(y)) = h(x,\pi_{F^\perp}(y))$, and from this fact we can conclude that $(\pi_{F^\perp}(\s))^\# = \s^\#\cap F^\perp$.
Proposition \ref{prop: new} then implies $\pi_{F^\perp}(\s)$ is an $\ord$-lattice in $F^\perp$ if and only if $\s^\#\cap F^\perp$ is an $\ord$-lattice in $F^\perp$. (Note that for the backwards direction we are using the fact that $\pi_{F^\perp}(\s)$ is an $\ord$-submodule of $(\s^\#\cap F^\perp)^\#$ satisfying $\R\pi_{F^\perp}(\s) = F^\perp$.)  Therefore by Lemma \ref{prop: reldet} (1), $\s\cap F$ is an $\ord$-lattice in $F$ if and only if $\s^\#\cap F^\perp$ is an $\ord$-lattice in $F^\perp$. 

Suppose now that $\s\cap F$ and $\s^\#\cap F^\perp$ are $\ord$-lattices in $F$ and $F^\perp$, respectively, and recall that $(\pi_{F^\perp}(\s))^\# = \s^\#\cap F^\perp$.  By Lemma \ref{prop: reldet} (2) we have the identity  
$$\det(\s) = \det(\s\cap{F})\det(\s^\#\cap F^\perp)^\#,$$
which by Proposition \ref{prop: new} can be rewritten as
$$\det(\s) = \det(\s\cap{F})\det(\s^\#\cap F^\perp)^{-1}{D_\ord}^{2(m-s)}.$$
\end{proof}

\section{Mordell's inequality for $\ord$-lattices}

Using the notation introduced in the previous section, define the $\ord$-Hermite constant $\gamma(\ord,rm)$ to be the supremum of the Hermite invariant 
$$\gamma(\s) = \frac{\mathrm{N}(\s)}{\det(\s)^{1/(rm)}}$$
for an $rm$-dimensional $\ord$-lattice $\s$.  
Observe that this constant is finite since the Hermite invariant of an $rm$-dimensional lattice is directly proportional to the $\left(\frac{rm}{2}\right)^\mathrm{th}$ power of its density, with the latter quantity bounded by $1$.  Moreover, by Mahler's compactness theorem (see \cite[p.~43]{Ma}) there exists an $rm$-dimensional $\ord$-lattice $\s$ with $\det(\s) = 1$ such that $\gamma(\ord, rm) = \gamma(\s)$, i.e., $\s$ is an optimal $rm$-dimensional $\ord$-lattice (note that here we are using the fact that the set of $rm$-dimensional $\ord$-lattices with determinant one is closed in the space of $rm$-dimensional lattices, and $\gamma(\ord, rm)$ is equal to the optimal value of norms of these lattices which must be bounded since $\gamma(\ord, rm)$ is finite).  

Below we use the results of Section 2 to prove an inequality relating the $\ord$-Hermite constants for dimensions $r(m-1)$ and $rm$, provided $m\geq 3$.  We choose to call this inequality a Mordell inequality for $\ord$-lattices due to its resemblance to the inequality $\gamma_{n-1}\leq \gamma_n^{(n-1)/(n-2)}$ in Mordell's theorem (see \cite[p.~41]{Ma}). In Mordell's theorem the constant $\gamma_n$ denotes Hermite's constant for dimension $n$, i.e., the optimal value of the Hermite invariant of an $n$-dimensional lattice.

\begin{theorem}\label{thm: mordell}
For each integer $m\geq 3$, 
\begin{equation}\label{mordell}
{\gamma{(\ord, rm)}}\leq{\gamma{(\ord,r(m-1))}}^{\frac{m-1}{m-2}}D_\ord^{\frac{1}{r(m-2)}}.
\end{equation}
Equality holds if and only if the $\ord$-dual lattice of every optimal $rm$-dimensional $\ord$-lattice $\s$ is also optimal, and for all minimal vectors $x\in{\s}^\#$ the relative $\ord$-lattices $\s\cap (\mathbb{K}x)^\perp$ and $\s^\#\cap \mathbb{K}x$ satisfy:

\begin{enumerate}
  \item $\mathrm{N}(\s \cap (\mathbb{K}x)^\perp) = \mathrm{N}(\s)$,
  \item $\s^\#\cap \mathbb{K}x = \ord x$,
  \item $\gamma(\s \cap (\mathbb{K}x)^\perp)=\gamma(\ord, r(m-1))$.
\end{enumerate}
\end{theorem}

\begin{proof}
Let $\s$ be an optimal $\ord$-lattice in $ E_{rm}$ with $\det(\s) = 1$,  
let $x$ be a minimal vector in the $\ord$-dual lattice $\s^{\#}$, and let $F$ denote the $(m-1)$-dimensional subspace ${(\mathbb{K} x)}^{\perp}$.  
By Proposition \ref{prop: intersection} the relative $\ord$-lattice $\s\cap F$ is an $\ord$-lattice in $F$ because $\s^\#\cap F^\perp$ is an $\ord$-lattice in $F^\perp = \mathbb{K}x$ (note that $\ord x\subseteq \s^\#\cap F^\perp$).  Then by the definition for the Hermite invariant of a lattice  
\begin{eqnarray*}
\gamma(\s)& = & \mathrm{N}(\s)\\
&\leq & \mathrm{N}(\s\cap F)\\
&= & \gamma(\s\cap F)\det(\s\cap F)^{\frac{1}{r(m-1)}}.
\end{eqnarray*}

We wish to bound the term $\det(\s\cap F)$ in the last equality by an expression involving the Hermite invariant of the $\ord$-dual lattice $\s^\#$.  To do this we first use the determinant identity in Proposition \ref{prop: intersection} with the definition of the Hermite invariant to get
\begin{eqnarray*}
\det(\s \cap F)& = &\det(\s)\det(\s^{\#} \cap F^\perp) {D_\ord}^{-2}\\
&\leq &\det(\ord x) {D_\ord}^{-2}.
\end{eqnarray*}
Then we let $\alpha_1,\dots,\alpha_r$ denote a $\Z$-basis for $\ord$ so that   
\begin{eqnarray*}
\det(\ord x)& = &\det((\langle\alpha_i x,\alpha_j x\rangle)_{1\leq i,j\leq r})\\
&= &\det\left(\left(\frac{1}{2}\left(h(\alpha_i x,\alpha_j x) + \overline{h(\alpha_i x,\alpha_j x)}\right)\right)_{1\leq i,j\leq r}\right)\\
&= &\det\left(\left(\frac{\mathrm{N}(x)}{2}\left(\alpha_i \overline{\alpha_j} + \alpha_j\overline{\alpha_i}\right)\right)_{1\leq i,j\leq r}\right)\\
&= &\mathrm{N}(\s^\#)^r{D_\ord}.
\end{eqnarray*} 
Observe that by Proposition \ref{prop: new} and the definition of $\gamma(\s^\#)$, the term $\mathrm{N}(\s^\#)^r$ in the last expression is equal to $\gamma(\s^\#)^r {D_\ord}^{2}$.  Hence we have shown that $\det(\ord x) = {\gamma(\s^\#)^r}{D_\ord}^3$, with which the upper bound for $\det(\s\cap F)$ computed above implies that    
$$\det(\s\cap F)\leq {\gamma(\s^\#)^r}D_\ord.$$  

Substituting the last upper bound for $\det(\s\cap F)$ into the initial expression involving $\gamma(\s)$ we obtain
\begin{eqnarray*}
\gamma(\s)& \leq &\gamma(\s\cap F)\left({\gamma(\s^\#)^r} {D_\ord}\right)^{\frac{1}{r(m-1)}}\\
&= &\gamma(\s\cap F) \gamma(\s^\#)^{\frac{1}{(m-1)}}{D_\ord}^{\frac{1}{r(m-1)}}\\
&\leq &\gamma(\ord, r(m-1)) \gamma(\ord, rm)^{\frac{1}{(m-1)}}{D_\ord}^{\frac{1}{r(m-1)}}.
\end{eqnarray*}
\ \\
Then, since $\gamma(\s) = \gamma(\ord,rm)$,
\begin{eqnarray*}
\gamma(\ord, rm)& \leq & \gamma(\ord, r(m-1)) \gamma(\ord, rm)^{\frac{1}{(m-1)}}{D_\ord}^{\frac{1}{r(m-1)}},
\end{eqnarray*}
and from this we obtain our desired inequality,
\begin{eqnarray*}\label{eqn: ineq}
\gamma(\ord, rm)& \leq & \gamma(\ord, r(m-1))^{\frac{m-1}{m-2}}{D_\ord}^{\frac{1}{r(m-2)}}.
\end{eqnarray*}

Observe that from the above proof of Inequality (3.1), equality holds if and only if all of the inequalities introduced in bounding $\gamma(\s)$ are tight, i.e., if and only if $\gamma(\s^\#) = \gamma(\ord, rm)$ and (1), (2) and (3) hold for every minimal vector $x\in\s^\#$ (the proof above uses an arbitrary minimal vector $x\in\s^\#$).   
More generally, this last statement holds for any optimal lattice $\s$ that does not necessarily have determinant $1$. 
\end{proof}

We choose to call the bound for $\gamma(\ord, rm)$ given by Inequality (3.1) the \it Mordell bound \normalfont for $rm$-dimensional $\ord$-lattices.  In Corollary \ref{cor: iterate} below we give an iterated version of this bound which  follows from a simple induction argument with Theorem \ref{thm: mordell} as the base case.

\begin{corollary}\label{cor: iterate}
For each pair of integers $s>m\geq2$,
\begin{equation*}\label{mordell}
{\gamma{(\ord, rs)}}\leq{\gamma{(\ord,rm)}}^{\frac{s-1}{m-1}}{D_\ord}^{\frac{s-m}{r(m-1)}}.
\end{equation*}

\end{corollary}

 \section{Eisenstein and Hurwitz Lattices}
Recall that Hurwitz lattices are lattices in quaternionic Hermitian space that are modules over the maximal $\Z$-order $\w = \Z[i,j,\frac{1+i+j+k}{2}]$ in $\left(\frac{-1,-1}{\Q}\right)$, and we say a $4m$-dimensional lattice (in Euclidean space) has a Hurwitz structure if it is isometric to a Hurwitz lattice.
Similarly, Eisenstein lattices are lattices in complex Hermitian space that are modules over the maximal $\Z$-order $\mathcal{E} = \Z[\frac{1+\sqrt{-3}}{2}]$ in $\Q(\sqrt{-3})$ (i.e., $\mathcal{E}$ is the ring of integers in $\Q(\sqrt{-3})$), and we say a $2m$-dimensional lattice has an Eisenstein structure if it is isometric to an Eisenstein lattice.

In this section we use the Mordell bounds for Eisenstein and Hurwitz lattices (Theorem \ref{thm: mordell}) to obtain upper bounds for the Hermite invariants of low-dimensional lattices with an Eisenstein or Hurwitz structure.  The upper bounds we obtain are compared in Tables 2 and 4 to the best upper bounds previously known and in several instances give an improvement. 
Note that for the dimensions in which the Eisenstein and Hurwitz-Hermite constants are not known, for comparison purposes we use the upper bounds for Hermite's constants in these dimensions proved by Henry Cohn and Noam Elkies in \cite{CE}; the author is not aware of any better bounds computed specifically for lattices with an Eisenstein or a Hurwitz structure.

\subsection{The Mordell Bound for Hurwitz Lattices} 
The densest known lattices with Hurwitz structures in dimensions $4m\leq 28$ are listed in Table 1 with their Hermite invariants\footnote{See \cite{NS} for further information on the lattices listed in Table 1.}.  The lattices listed in this table for dimensions $4$, $8$, $12$, and $24$ have all been proven optimal as lattices with a Hurwitz structure, and with the exception of dimension $12$ this is a corollary to their proven optimality as lattices; the $12$-dimensional lattices listed in Table 1 have been proven optimal as lattices with Hurwitz structures by Fran\c cois Sigrist and David-Olivier Jaquet-Chiffelle, and a summary of their calculations is given in \cite{Si}.\footnote{The details of Sigrist and Jaquet-Chiffelle's calculations can also be found in Achill Sch$\mathrm{\ddot{u}}$rmann's paper \cite{Sch}.} 

\begin{table}\label{t3}
\begin{center}
\caption{The densest known lattices with Hurwitz structures in dimensions $4m\leq 28$}
\begin{tabular}{c||c|r}\hline
$4m$ & Lattice(s) & Hermite Invariant \\
\hline 

$4$ & $\s_4 = D_4$ & $\sqrt{2} \approx 1.41421$\\
$8$ & $\s_8 = E_8$  & $2$ \\
$12$ & $\s_{12}^{min}$, $\s_{12}^{max}$  & $2^{7/6}\approx 2.24492$  \\
$16$ & $\s_{16}$ & $2^{3/2} \approx 2.82843$ \\
$20$ & $\s_{20}$ & $2^{17/10} \approx 3.24901$ \\
$24$ & $\s_{24}$ & $4$ \\
$28$ & $\s_{28}$, $LL_{28}$ & $2^{27/14} \approx 3.80678$\\
\hline
\end{tabular}
\ \\
\end{center}
\end{table}

In Table 2 we list the Mordell bounds for Hurwitz lattices in dimensions $4m\leq 28$ alongside the best upper bound for $\gamma(\w,4m)$ previously known.  Note that each of the Mordell bounds listed have been computed iteratively; i.e.,  if the Mordell bound for $4(m-1)$-dimensional Hurwitz lattices is an improvement on the best known upper bound previously known for $\gamma(\mathcal{H},4(m-1))$, then it is used to compute the Mordell bound for $4m$-dimensional Hurwitz lattices.  We have also included in Table 2 a conjectured Mordell bound for $24$-dimensional Hurwitz lattices which is obtained by replacing $\gamma(\w, 20)$ with $\gamma(\s_{20})$ in Inequality (\ref{thm: mordell}).  

\begin{table}\label{t4}
\begin{center}
\caption{The Mordell bound and conjectured Mordell bound for the Hurwitz Hermite constant for dimensions $4m\leq 28$}
\begin{tabular}{c||c|r|r}\hline
$4m$ & Best Known & Mordell Bound & Conjectured  \\
\ &Upper Bound & & Mordell Bound \\
\hline 
$4$ & $\sqrt{2}\approx 1.41421$ & -- & -- \\
$8$ & $2$  & -- & -- \\
$12$ & $2^{7/6}\approx 2.24492$ & $2^{3/2} \approx 2.82843$ & -- \\
$16$ & $3.02639$ & $\mathbf{2^{3/2} \approx 2.82843}$ & -- \\
$20$ & $3.52006$ & $2^{11/6} \approx 3.56359$ & --\\
$24$ & $4$ & $4.21390$ & $4$ \\
$28$ &  $4.48863$   & $2^{23/10} \approx 4.92458$ & -- \\
\hline
\end{tabular}
\ \\
\end{center}
\end{table}

Notice that the Mordell bound for $16$-dimensional Hurwitz lattices is not only an improvement on the best upper bound for $\gamma(\w, 16)$ previously known, but it is equal to the Hermite invariant of the $16$-dimensional Barnes-Wall lattice $\s_{16}$. Therefore we have proved the following theorem.

\begin{theorem}\label{thm: BW}
The $16$-dimensional Barnes wall lattice $\s_{16}$ has optimal density as a $16$-dimensional lattice with a Hurwitz structure.  Moreover, every optimal $16$-dimensional lattice with a Hurwitz structure contains either $\s_{12}^\mathrm{min}$ or $\s_{12}^\mathrm{max}$ as a $12$-dimensional section with identical norm.
\end{theorem}

Note that we have to be careful when translating the conditions in Theorem \ref{thm: mordell} to $4m$-dimensional lattices with a Hurwitz structure because it may be possible for a $4m$-dimensional lattice to have two or more inequivalent Hurwitz structures; that is, a lattice may be isometric to two different Hurwitz lattices such that there does not exist an isometry between them that is an $\w$-module isomorphism.  For example it is currently an open problem to determine if $\s_{16}$ has two inequivalent Hurwitz structures (see \cite[p.~279]{Ma}).  With the Hurwitz structure presently known for $\s_{16}$, all of the $12$-dimensional sections perpendicular to a minimal vector in $\s_{16}^\#$ are $\s_{12}^\mathrm{max}$.
Jacques Martinet has shown in \cite[Ch.~8]{Ma} that if another inequivalent Hurwitz structure exists for $\s_{16}$, then all of the $12$-dimensional sections perpendicular to a minimal vector in $\s_{16}^\#$ (computed with respect to the new structure) are $\s_{12}^\mathrm{min}$.

One final observation we wish to make concerning Table 2 is that the conjectured Mordell bound for $24$-dimensional Hurwitz lattices is tight, and from this we can use Theorem \ref{thm: mordell} to make conjectures about the conditions satisfied by certain $20$-dimensional sections of the $24$-dimensional Leech lattice as was done for the $16$-dimensional Barnes-Wall lattice above.  However, unlike the latter case, the Leech lattice has been proven to have a unique Hurwitz structure by H.-G. Quebbemann (see \cite{Qu}).

\subsection{The Mordell Bound for Eisenstein lattices}
We now list in Table 3 the densest known lattices with Eisenstein structures in dimensions $2m\leq 26$ with their Hermite invariants\footnote{See \cite{NS} for further information on the lattices in Table 3.}$^,$\footnote{All lattices listed in Table 3, with the exception of the $26$-dimensional laminated Eisenstein lattice $\s_{26}$, are unique as lattices and have a unique Eisenstein structure; see \cite{CS2} and \cite{CS3} for more information on the non-uniqueness of $\s_{26}$ as a lattice and for the uniqueness of the Eisenstein structures on the other lattices.}.  All of the lattices in Table 3 are among the densest known lattices in their dimension; however, we note that this table is not a complete list of the densest known lattices for the dimensions shown because the density of the Conway-Borcherds lattice $T_{26}$ is equal to that of $\s_{26}$ and the former lattice does not have an Eisenstein structure.  The lattices listed for dimensions $2m\leq 8$ and for dimension $24$ have all been proven optimal as ordinary lattices (and hence are optimal as lattices with an Eisenstein structure), and the lattice $\s_{10}$ has recently been proven optimal as a lattice with an Eisenstein structure by Achill Sch$\ddot{\mathrm{u}}$rmann; see \cite{Sch} for details.

\begin{table}\label{t1}
\begin{center}
\caption{The densest known lattices with Eisenstein structures in dimensions $2m\leq 26$}
\begin{tabular}{c||c|r}\hline\\
$2m$ & Lattice(s)
& Hermite Invariant \\
\hline 
$2$ & $\s_2 = A_2$ & $2/\sqrt{3} \approx 1.15470$\\
$4$ & $\s_4 = D_4$ & $\sqrt{2} \approx 1.41421$\\
$6$ & $\s_6 = E_6$ & $2/3^{1/6} \approx 1.66537$ \\
$8$ & $\s_8 = E_8$  & $2$ \\
$10$ & $\s_{10}$ & $2^{6/5}/3^{1/10} \approx 2.05837$ \\
$12$ & $K_{12}$  & $4/\sqrt{3} \approx 2.30940$ \\
$14$ & $\s_{14}$ & $2^{10/7}/3^{1/14} \approx 2.48864$ \\
$16$ & $\s_{16}$ & $2^{3/2} \approx 2.82843$ \\
$18$ & $\s_{18}$ & $2^{5/3}/3^{1/18} \approx 2.98683$ \\
$20$ & $\s_{20}$ & $2^{17/10} \approx 3.24901$ \\
$22$ & $\s_{22}$ & $2^{21/11}/3^{1/22} \approx 3.57278$ \\
$24$ & $\s_{24}$ & $4$ \\
$26$ & $\s_{26}$ & $4/3^{1/26} \approx 3.83450$ \\
\hline
\end{tabular}
\ \\
\end{center}
\end{table}

In Table 4 we list the Mordell bounds for Eisenstein lattices in dimensions $2m\leq 26$ next to the best upper bound previously known for $\gamma(\mathcal{E}, 2m)$.  Note that the Mordell bounds listed in this table are computed iteratively, as was done for the Hurwitz case.  Included in Table 4 is a column listing the conjectured Mordell bound for Eisenstein lattices in the dimensions for which the optimal density of an Eisenstein lattice with one less complex dimension is not known.  These conjectured bounds are computed using Inequality (\ref{thm: mordell}) in Theorem \ref{thm: mordell} with the Hermite invariants of the densest known $2(m-1)$-dimensional Eisenstein lattices replacing $\gamma(\mathcal{E}, 2(m-1))$.

\begin{table}\label{t2}
\begin{center}
\caption{The Mordell bound and conjectured Mordell bound for the Eisenstein Hermite constant for dimensions $2m\leq 26$}
\begin{tabular}{c||r|r|r}\hline
Dimension & Best Known  & Mordell Bound & Conjectured  \\
 $2m$     & Upper Bound &               & Mordell Bound \\
\hline 
$2$  & $\sqrt{3}/2\approx 1.15470$ & --                         & -- \\
$4$  & $\sqrt{2}\approx 1.41421$ & --                         & -- \\
$6$  & $2/3^{1/6}\approx 1.66537$ & $\sqrt{3} \approx 1.73205$ & -- \\
$8$  & $2$       & $2$                        & -- \\
$10$ & $2.05837$ & $2\cdot 3^{1/6} \approx 2.40187$ & -- \\
$12$ & $2.52179$ & $\mathbf{2^{5/4} \approx 2.37841}$ & -- \\
$14$ & $2.77580$ & $\mathbf{2^{13/10}\cdot3^{1/10} \approx 2.74822}$ & $\mathbf{2^{11/5}/\sqrt{3} \approx 2.65281}$  \\
$16$ & $3.02639$ & $3.17552$ & $\mathbf{2^{3/2} \approx 2.82843}$ \\
$18$ & $3.27433$ & $3.47300$ & $\mathbf{2^{11/7}\cdot 3^{1/14} \approx 3.21460}$  \\
$20$ & $3.52006$ & $3.72996$ & $\mathbf{2^{7/4} \approx 3.36359}$  \\
$22$ & $3.76404$ & $3.98416$ & $\mathbf{2^{16/9}\cdot3^{1/18} \approx 3.64478}$ \\
$24$ & $4$ & $4.23616$ & $4$ \\
$26$ & 4.24804   & $2^{23/11}\cdot 3^{1/22} \approx 4.47831$ & -- \\
\hline
\end{tabular}
\ \\
\end{center}
\end{table}

The bold entries in Table 4 indicate an improvement (actual or conjectured) on the best known upper bound for $\gamma(\mathcal{E},2m)$. Specifically, we have improved on the best known upper bounds for the Hermite invariants of Eisenstein lattices in dimensions $12$ and $14$, with the Mordell bound computed for $\gamma(\mathcal{E},12)$ being much closer to the Hermite invariant of the Coxeter-Todd lattice $K_{12}$ than the upper bound from \cite{CE}. 
Unfortunately with the improved bound for $\gamma(\mathcal{E},12)$ we are unable to conclude that the Coxeter Todd lattice $K_{12}$ is optimal as a $12$-dimensional lattice with an Eisenstein structure, despite the fact that it is widely believed to be optimal as a $12$-dimensional lattice.  However, one could have anticipated that this bound would not be sharp because $K_{12}$ does not contain the lattice $\s_{10}$ as a $10$-dimensional section, but rather contains the lattice $K_{10}'$.  For if the Mordell bound for $\gamma(\mathcal{E},12)$ were tight, then by Theorem \ref{thm: mordell} any optimal lattice would contain $\s_{10}$ as a $10$-dimensional section because $\s_{10}$ is the unique optimal $10$-dimensional lattice with an Eisenstein structure.

Moving on to the conjectured Mordell bounds for Eisenstein lattices displayed in Table 4, note that if one could prove $\s_{14}$ is optimal as a $14$-dimensional lattice with an Eisenstein structure, then the Mordell bound for $16$-dimensional Eisenstein lattices would imply that the $16$-dimensional Barnes-Wall lattice $\s_{16}$ is optimal as a lattice with an Eisenstein structure.
Moreover, due to the conjectured tightness of the Mordell bounds for $16$ and $24$-dimensional Eisenstein lattices,   we can make conjectures about optimal $16$ and $24$-dimensional Eisenstein lattices, as was done in the Hurwitz case above for these same dimensions.

\section{Remarks and Open Problems}

Even though the focus of this paper is on lattices with Eisenstein and Hurwitz structures, we wish to emphasize that Theorem \ref{thm: mordell} applies to other types of $\ord$-lattices.  Other $\ord$-lattices one can consider include the $2m$-dimensional Gaussian lattices ($\ord = \Z[i]$) and the $4m$-dimensional $\mathcal{J}$-lattices where $\mathcal{J}$ is the subring $\Z[1,i,\frac{1+\sqrt{3}j}{2},\frac{i+\sqrt{3}k}{2}]$ in $\h$ corresponding to the maximal $\Z$-order $\Z[1,i,\frac{1+j}{2},\frac{i+k}{2}]$ in the rational positive definite quaternion algebra $$\left(\frac{-1, -3}{\Q}\right) = \lbrace a+bi+cj+dk: a,b,c,d\in\Q \text{\ and\ } i^2 = -1, j^2 = -3, ij = -ji = k\rbrace.$$ 
Many of the densest known even-dimensional lattices have a Gaussian structure, and the $4m$-dimensional $\mathcal{J}$-lattices are of interest because they include the  Coxeter-Todd lattice $K_{12}$, the densest known 12-dimensional lattice, which curiously does not have the structure of a Hurwitz lattice; see \cite{Gr} concerning the existence of a $\mathcal{J}$-lattice structure on $K_{12}$.  Note that the $\mathcal{J}$-Hermite constant for dimension $8$ is equal to $\gamma(E_{8})$; unfortunately $K_{12}$ does not contain $E_8$ as an $8$-dimensional section (see Proposition 8.7.9 in \cite{Ma}), and so we cannot use Theorem \ref{thm: mordell} to conclude that $K_{12}$ is optimal as a $12$-dimensional $\mathcal{J}$-lattice.

It may be possible to use Theorem $3.1$ to prove that an $r(m-1)$-dimensional lattice \it does not \normalfont have a particular $\ord$-lattice structure.  As pointed out to the author by Jacques Martinet, if we know the value of $\gamma(\ord, rm)$ but do not know the value of $\gamma(\ord, r(m-1))$, then we can replace $\gamma(\ord, r(m-1))$ in the inequality by the Hermite invariant of an $r(m-1)$-dimensional lattice $\s$ and check if the inequality remains valid.  If the inequality is no longer valid, then this implies that $\s$ does not have an $\ord$-lattice structure.  It would be interesting to find examples of lattices that can be proven not to have certain $\ord$-lattice structures using this method.  

By Theorem 8.7.2 in \cite{Ma}, the $8$-dimensional lattice $E_8$ has an $\ord$-lattice structure over every maximal order $\ord$ in an imaginary quadratic number field.  Given that the Mordell inequality for Eisenstein and Gaussian lattices in dimension $8$ is sharp (note that the densest $6$-dimensional Gaussian lattice is $D_6$; see \cite{Sch}), it is natural to wonder if the Mordell inequality is sharp for $\ord$-lattices in dimension $8$ whenever $\ord$ is a maximal order in an imaginary quadratic number field.  Even if this is not true in general, it still would be interesting if one could characterize the maximal orders $\ord$ for which this inequality is sharp.

\section{Acknowledgments}
First and foremost the author thanks her advisor, Henry Cohn, for introducing her to this area of mathematics and for his constant support of her research pursuits.  Next the author thanks Achill Sch$\ddot{\mathrm{u}}$rmann for sharing his results on perfect forms and for making himself available to answer questions and to offer suggestions on how this paper could be improved.  The author also thanks the Hausdorff Research Institute for Mathematics for hosting the workshop titled \lq\lq Experimentation with, Construction of, and Enumeration of Optimal Geometric Structures\rq\rq.  This workshop provided the opportunity to present a preliminary version of this work and to get many helpful comments and suggestions on how it could be generalized.  In particular the author thanks workshop participants Jacques Martinet, Gabriele Nebe, and Renaud Coulangeon for their interest in her research and for suggesting that she use localization to generalize her results to include $\ord$-lattices where the maximal $\Z$-order $\ord$ is not a left principal ideal domain and hence every $\ord$-lattice is not necessarily a free $\ord$-module.  Additionally, the author thanks Jacques Martinet for his later suggestions and comments and for providing references on maximal orders.

\bibliographystyle{amsalpha}

\end{document}